\documentclass[11pt, a4paper]{article}

\usepackage{tikz} 

\usepackage[centertags]{amsmath} 
\usepackage{amsthm,amssymb,amscd,amsfonts}
\usepackage{mathtools}
\usepackage[all]{xy}
\usepackage{enumitem}

\usepackage[english]{babel}
\usepackage[bookmarks,breaklinks]{hyperref}

\usepackage{breqn}

\usepackage{mathrsfs}  

\theoremstyle{plain}

\newtheorem{thm}{Theorem}[section]
\newtheorem{cor}[thm]{Corollary}
 
\newtheorem{lemma}[thm]{Lemma}
\newtheorem{prop}[thm]{Proposition} 
\newtheorem{defn}[thm]{Definition}
\newtheorem{assumption}{Assumption}

\theoremstyle{remark}
\newtheorem*{rem}{Remark}

\newcommand{\Oo}{\mathcal{O}}

\newcommand{\R}{\mathbb{R}} 
\newcommand{\m}{\mathfrak{m}}

\newcommand{\Z}{\mathbb{Z}}
 
\newcommand{\C}{\mathbb{C}}

\newcommand{\an}{\operatorname{an}}
\newcommand{\Cmer}{\operatorname{\mathbb{C}_{mer}}}
\renewcommand{\div}{\operatorname{div}}
\newcommand{\ord}{\operatorname{ord}}
\newcommand{\red}{\operatorname{red}}
\newcommand{\sing}{\operatorname{sing}}
\newcommand{\val}{\operatorname{val}}

\newcommand{\diam}{\operatorname{diam}}
\newcommand{\dis}{\operatorname{dis}}

\newcommand{\St}{\operatorname{St}}
\newcommand{\Sk}{\operatorname{Sk}}
\newcommand{\eps}{\varepsilon}

\newcommand{\wo}{\setminus}

\newcommand{\suchthat}[2]{\{\ #1\ \mid\ #2\ \}}

\newcommand{\wt}{\operatorname{wt}}

\newcommand{\Arg}{\operatorname{Arg}}

\newcommand{\Xx}{\mathscr{X}}
\newcommand{\Yy}{\mathscr{Y}}
\newcommand{\XX}{\mathcal{X}}

\newcommand{\D}{\mathbb{D}}
\newcommand{\what}[1]{\widehat{#1}}

\begin{document}

\bibliographystyle{alpha}

\begin{center}
  {\Large\bf Gromov-Hausdorff limits of flat Riemannian surfaces}\\[2ex]
\end{center}

\begin{center}
  {\normalsize Dmitry Sustretov \footnote{ This project has received
      funding from the European Union's Horizon 2020 research and
      innovation programme under the Marie Sk\l{}odowska-Curie grant
      agreement No. 843100 (NALIMDIF). The author was supported by the
      European Research Council under the European Community's Seventh
      Framework Programme (FP7/2007-2013) with ERC Grant Agreement
      nr. 615722 MOTMELSUM, and by Labex CEMPI
      (ANR-11-LABX-0007-01). 
  }
  {\small \hspace{0.15\linewidth}}}
\begin{minipage}[t]{0.85\linewidth}
  \begin{center} {\bf Abstract}
  \end{center}
  I study Gromov-Hausdorff limits of complex curves endowed with
  singular flat metrics of constant diameter. I formulate a criterion
  that the limit is collapsed in terms of a certain piecewise affine
  weight function on the dual intersection complex of a semi-stable
  model of the degeneration introduced by Kontsevich and Soibelman. I
  describe the collapsed and non-collapsed limits, which are,
  respectively, metric graphs and finite collections of complex curves
  with flat metrics glued along finitely many points. I show that the
  collapsed limit of any positive genus can occur.
\end{minipage}

\end{center}

\tableofcontents


\section{Introduction}

In the paper \cite{ks-affine} Kontsevich and Soibelman formulate a
series of conjectures about the shape of Gromov-Hausdorff limits of
certain families of complex manifolds endowed with Ricci flat
metrics. These conjectures are motivated by mirror symmetry and in
particular by the authors' approach to the SYZ conjecture.  One
considers germs of holomorphic families of compact Calabi-Yau
manifolds parametrized by points of a punctured disc having maximally
unipotent action of the monodromy on the middle cohomology and with a
relatively ample line bundle on the total space of the family. For
each element of the family one picks the Ricci flat metric with the
K\"ahler class equal to the first Chern class of the polarizing line
bundle and normalized so that the diameter is
constant. Gromov-Hausdorff limits of such families are then
conjectured by Kontsevich and Soibelman to carry a singular affine
manifold structure with respect to which the limit metric satisfies
the real Monge-Ampere equation. The real dimension of the limit
manifold is half the real dimension of the elements of the family, so
we speak about \emph{collapsed} limits.

Alternatively, the limit manifold together with the singular affine
structure can be recovered from the non-archimedean analytification
$X^{\an}$ (in the sense of Berkovich) of the variety $X$ over the
non-archimedean field $\Cmer$ of germs of complex functions
meromorphic at 0 (\cite[\S 5]{ks-affine}). As a topological space it
is a closed subset of $X^{\an}$, the minimality locus of a certain
weight function associated to a canonical form; this closed subset is
called the essential skeleton of $X$, denoted $\Sk(X)$. On a variety
with trivial canonical bundle weight functions associated to different
canonical forms differ by a constant and so the minimality locus does
not depend on this choice. As it turns out, there exists a (generally
speaking, not canonical) retraction $X^{\an} \to \Sk(X)$ \cite{nx-ess}
which is a fibration over an open dense subset of $\Sk(X)$, with the
fibre isomorphic to a non-archimedean torus \cite{nx-yu}; this endows
$\Sk(X)$ with an integral affine structure away from the
discriminant locus. Kontsevich and Soibelman conjecture
(\cite[Conjecture~3, \S 5]{ks-affine}) that $\Sk(X)$ is isomorphic as
a mainfold with affine structure with singularities to the
Gromov-Hausdorff limit described in the previous paragraph.

Collapsed Gromov-Hausdorff limits of Ricci-flat hyperk\"ahler
manifolds have been extensively studied in a slightly different
setting: one fixes a holomorhic fibration of a single hyperk\"ahler
manifold, with a generic fibre Abelian variety, and considers a
variation of the K\"ahler class, making it tend to the boundary of the
K\"ahler cone \cite{gross2000large}, \cite{gtz-collapsing},
\cite{gtz-collapsing2}, \cite{tz-hk}.  Boucksom and Jonsson
\cite{boucksom2016tropical} show that volume forms on the fibres of
the degenerating family converge in a certain sense to a measure
supported on $\Sk(X)$. In the paper \cite{odaka-oshima} Odaka
and Oshima describe Gromov-Hausdorff limits of K3 surfaces in the set
up of \cite[Conjecture~3, \S 5]{ks-affine} (though they do not discuss
the relationship with $\Sk(X)$). Finally, I would like to
mention a recent paper \cite{dhl} of Ducros, Hrushovski and Loeser
where they propose a framework for asymptotic integration that is
motivated by the said Conjecture.

In this paper we consider the Gromov-Hausdorff limits of families of
complex curves of genus $\geq 1$ endowed with flat pseudo-K\"ahler
metrics. If one compares to the set up of Kontsevich and Soibelman, I
relax the assumption on the triviality of the tangent bundle and allow
the metric to have conical singularities in finitely many points. We
rescale the metrics with the K\"ahler form
$\dfrac{i}{2}(\Omega_s \wedge \bar\Omega_s)$ (where $\Omega$ is a
given relative 1-form) on the fibres $X_s$ of the family $X$ so that
$\diam X_s \equiv 1$.

We show that there are two possibilites for the limit. In the
collapsed case the limit is a metric graph which can be canonically
represented as a certain quotient of the dual intersection complex of
the special fibre of a semi-stable model of the degeneration. The
quotient is defined in terms of minimality locus of the weight
function associated to the form $\Omega$, defined in \cite{ks-affine}
and further studied in \cite{mn-weight, nx-ess} and \cite{temkin}. In
the non-collapsed case, the limit is a union of flat surfaces glued
along finitely many points at which the metric is singular, the gluing
is determined by the minimality locus of the weight function.

The main results of this paper are Theorems~\ref{limit-collapsed} and
\ref{limit-non-collapsed} which describe the collapsed and
non-collapse d limits. In Proposition~\ref{collapsed-any-genus} a
series of degenerations of curves of genus $2k+1$ which give rise to
collapsed limits, metric graphs of any genus $k \geq 1$, are
constructed.  The technical heart of the paper is Section~\ref{limits}
where a neighbourhood of the special fibre of a model of the
degeneration is covered by charts of a special form and estimates on
the lengths of shortest geodiscs are derived. Section~\ref{prelim}
provides background information about dual intersection complexes and
the weight function.

\textbf{Acknowledgements}. I would like to thank Ehud Hrushovski for
telling me about the conjectures of Kontsevich and Soibelman and
suggesting the idea to use the Robinson field for the description of
the Gromov-Hausdorff limit which was crucial in the initial version of
this work. I would like to thank Grisha Papayanov, J\'er\^ome Poineau,
Matthew Stevenson, Ilya Tyomkin, Misha Verbitsky and Eric Walsberg for
helpful discussions.  I also would like to thank I.H.\'E.S. and Max
Planck Institute, where the work on this paper was carried out, for
perfect working conditions.

\section{Background}
\label{prelim}

\subsection{Gromov-Hausdorff distance}

We refer the reader to \cite[Chapter~7]{bbi} for a comprehensive 
introduction into the Gromov-Hausdorff distance. Let $X$ be a metric
space with the metric $d$ and let $A, B \subset X$ be two subsets; the
\emph{Hausdorff distance} between $A$ and $B$ is the infimum of real
numbers $\eps > 0$ such that $B_{\eps}(X) \subset Y$ and
$B_{\eps}(Y) \subset X$, where $B_{\eps}(A)$
$$
B_{\eps}(A) = \suchthat{x \in X}{\exists a \in A\ \ d(x,a) < \eps}
$$
denotes the $\eps$-neighbourhood of a set $A$. The
\emph{Gromov-Hausdorff distance} between two metric spaces $(X,d)$ and
$(Y,d')$ is the infimum of Hausdorff distances between $X$ and $Y$
over all metric spaces $Z$ and all possible isometric embeddings of
$X \hookrightarrow Z$, $Y \hookrightarrow Z$. Note that finite metric
spaces are dense in the space of (isometry classes of) compact metric
spaces with the Gromov-Hausdorff metric.

If $X$ is a complex curve and $\Omega$ is a holomorphic 1-form on $X$
then $\omega=\dfrac{i}{2}(\Omega \wedge \bar\Omega)$ defines a
positive semi-definite (1,1)-form an $X$, so if $I$ the complex
structure then $g(x,y)=\omega(x, Iy)$ is a pseudo-Riemannian
metric. Such pseudo-Riemannian surfaces are locally isometric to the
Euclidean plane away from the zeroes of $\Omega$, where they have
conical singularities, and have trivial holonomy. They are also called
\emph{translation surfaces} since they can be glued from polygonal
domains on a plane via identification of opposite sides by
translations (see, for example, \cite{zorich2006flat}). The shortest
geodesic metric on such spaces is a complete inner metric.

Let $X$ be a variety over the field $\Cmer$ of germs of functions
meromorphic at 0. To give $X$ is the same as to give a germ of a
family $\XX \to \D_\eps^\circ$ over the punctured disc
$\D_\eps^\circ=\suchthat{x \in \C}{0 < |x| < \eps}$ for some
sufficiently small $\eps$. Assume that the genus $g(X)$ of $X$ is
greater than or equal to 1, and let
$\Omega \in H^0(X, \Omega_{X/\Cmer})$ be given. Denote
$\omega_s=\dfrac{i}{2}(\Omega_s \wedge \bar\Omega_s)$ the K\"ahler
forms on each fibre $X_s$ for $0 < |s| < \eps$ and let
$\tilde{\omega}_s=\omega_s/\diam(\Xx_s, \omega_s)^2$ be the
rescaled K\"ahler form (so that
$\diam (\Xx_s, \tilde{\omega}_s) \equiv 1$). Drawing analogy with the
conjecture of Kontsevich and Sobelman one can ask:
 
 \textbf{Question}: what is the limit in the Gromov-Haudorff metric of
 $(X_s, \tilde{\omega}_s)$ as $s \to 0$?

 The answer will be given in the Section~\ref{shape}.

\subsection{Dual intersection complexes and the weight function}
\label{dual}

Everything in this section is valid over a discretely valued field $K$
with a value ring $R$ ($K=\Cmer$ in the rest of the article), and the
valuation $v_K: K^\times \to \R$. For a scheme $\Xx$ over $R$ we
denote $\Xx_0$ the scheme-theoretic fibre over the closed point of
$R$.

Recall that a \emph{model} $\Xx$ is a flat scheme over $R$ such that
$\Xx \otimes_R K \cong X$. It is called an \emph{snc model} if
$(\Xx_0)_{\red}$ is a divisor with strict normal crossings.

\begin{defn}[Dual interesection complex]
  Let $\Xx$ be an snc model of a projective curve $X$, and let the
  central fibre $\Xx_0$ equal $\sum_{i=1}^m N_i E_i$ where
  $E_1, \ldots, E_m$ are irreducible components of $\Xx_0$. The
  \emph{dual intersection complex $\Delta_\Xx$} of the special fibre
  $\Xx_0$ is a metric graph that has vertices $[E_1], \ldots, [E_m]$
  and edges $[\sigma]$ of length $l(\sigma)=(N_i,N_j)/N_i N_j$ for
  each point $\sigma \in E_i \cap E_j$.
\end{defn}

For any edge $\sigma$ we denote $\partial \sigma$ the set of its ends.
For any vertex $[E_i] \in \Delta_\Xx$ we will denote $\St([E_i])$ the
\emph{star of $[E_i]$}, the set of edges $\sigma$ such that
$[E_i] \in \partial \sigma$.

The dual intersection complex of any snc model embeds into the
\emph{Berkovich analytification $X^{\an}$} of $X$. Let $v_K$ be the
valuation on the base field $K$. As a topological space the
analytification $X^{\an}$ is defined to be the set of pairs
$$
X^{\an} := \suchthat{(\xi, v)}{\xi \in X, v: K(\xi) \to \R
  \textrm{ valuation }, v|_K = v_K}
$$
with the weakest topology that makes evaluation maps $v \mapsto v(f)$
continuous for any $f \in K[U], U \ni x$.

The construction of the embedding
$\Delta(\Xx) \hookrightarrow X^{\an}$ goes back to \cite{berko}, we
recall here a more direct approach following
\cite[Proposition~3.1.4]{mn-weight}, \cite[Section~3]{bfj-siminag}. We
identify a face $\sigma$ joining two components $E_i$ and $E_j$ with
an interval in $\R_{\geq 0}^2$ given by the equation
$N_j x + N_i y = 1$ (note that the metric on $\sigma$ has nothing to
do with the Euclidean metric on this interval). A point with
coordinates $(\alpha, \beta)$ is identified with a quasi-monomial
valuation as follows. Let $x,y \in \what{\Oo_{\Xx, \sigma}}$ be some
local parameters at an intersection point $\sigma \in E_i \cap E_j$
such that $x$ is a local equation for $E_i$ and $y$ is a local
equation for $E_j$.  Define the following valuation:
$$
v_\alpha: K(X)^\times \to \R\qquad f \mapsto \min_{f_{ij} \neq 0}
\alpha i + \beta j
$$
where $f_{ij}$ are the coefficients of an expansion
$f = \sum_{i,j=1}^\infty f_{ij} x^i y^j$; by
\cite[Proposition~2.4.6]{mn-weight} this definition does not depend on
the choice of the local equations for $E_i, E_j$. One can show
(\cite[Proposition~3.2.2 ]{mn-weight}) that if $f$ has no zeroes or
poles that contain $\sigma$ then $v_\alpha(f)$ is an affine function
of $\alpha$.

\begin{rem}
  The part of the topological space $X^{\an}$ that consists of
  valuations on the function field of $X$ can be metrized, see, for
  example, \cite[5.58]{bpr}. The image of the embedding
  $\Delta_{\Xx} \to X^{\an}$ clealy lies in this part, and one can
  show that the embedding is isometric with respect to this metric.
\end{rem}

Let $\Yy$ be the blow-up of a point $\sigma \in E_i \cap E_j$. Then
$\Delta_\Yy$ is obtained from $\Delta_\Xx$ by the subdivision of the
edge that joins $[E_i]$ and $[E_j]$ and that corresponds to the
intersection of $E_i$ and $E_j$ that has been blown up; the new point
is the divisorial valuation correspnding to the exceptional divisor of
the blow-up. If $\Yy$ is a blow-up of a smooth point $x \in E_i$ of
$\Xx_s$ then $\Delta_\Yy$ is obtained by adjoining an edge to $[E_i]$
in $\Delta_\Xx$, so the latter dual intersection complex can be
naturally regarded as a subgraph of $\Delta_\Yy$.

For a blow-up $f: \Yy \to \Xx$ with the exceptional divisor $E$ there
exists a natural map $r: \Delta_\Yy \to \Delta_\Xx$ which retracts the
edge containing the point $[E]$ if $f(E) \not\in (\Xx_0)_{\sing}$, 
or sends it to the barycenter of the interval joining $[E_i]$ and
$[E_j]$ if $f(E) \subset E_i \cap E_j$. Since any model $\Yy$ that
dominates $\Xx$ is obtained as a sequence of blow-ups of points in the
central fibre, the map $r$ can be defined as a composition of such
maps for any dominant $\Yy \to \Xx$.  The retraction map can be
defined in a less ad hoc and still explicit way in any dimension, see
\cite[Theorem~3.1]{bfj-siminag}, \cite[Proposition~3.1.7]{mn-weight}.

Let $\Omega \in H^0(X, \Omega^1_X)$. Define the weight function
$\wt_\Omega: \Delta_\Xx \to \R \cup \{+\infty\}$ as the function that
takes the following values on the divisorial valuations
$$
\wt_{\Omega}([E_i]) = \dfrac{1 + \ord_{E_i}(\Omega)}{N_i}
$$
where $\ord_{E_i}(\Omega)$ is the order of vanishing at the divisor
$E_i$ of $\Omega$ regarded as the rational section of the relative
canonical bundle $\omega_{\Yy/R}$ of a model that has $E_i$ as one of
the components of the central fibre. By Proposition~4.2.4
\cite{mn-weight} $\wt_\Omega$ is well-defined (i.e. does not depend on
the model $\Yy$) and by Proposition~4.4.5 \emph{loc.cit} its extension
by continuity to the whole of $X^{an}$ gives rise to a function that
is piece-wise affine on the faces of $\Delta_\Xx$ for any snc model
$\Xx$.

It follows from this definition that the weight function is
compatible with the embeddings of dual intersection complexes: if
$f: \Yy \to \Xx$ is a dominant morphism of snc models, then
$$
\wt_\Omega^\Yy|_{\Delta_\Xx} = \wt_\Omega^\Xx
$$
Furthermore, by \cite[Proposition~3.1.6]{mn-weight} if
$r: \Delta_\Yy \to \Delta_\Xx$ is the retraction then
$\wt_\Omega^\Yy(v) \geq \wt_\Omega^\Yy(r(v))$ .

An alternative treatment of the weight function using the
non-archimedean analytic techniques can be found in \cite{temkin}.

\section{Limits of flat curves}
\label{limits}

\subsection{Models and charts}
\label{models-charts}

From now on $K=\Cmer$.

\begin{defn}[Snc model of a pair]
  A pair $(\Xx, \Omega')$ of a model $\Xx$ and a relative
  1-form $\Omega'$ such that $\Omega'|_X = \Omega$ is called
  a \emph{model of the pair $(X,\Omega)$}. It is called an \emph{snc
    model of the pair $(X,\Omega)$} if the reduction of
  $\div(\Omega') \cup \Xx_s$ is snc.
\end{defn}

\begin{lemma}
  \label{make-snc}
  Let $(\Xx,\Omega')$ be a model of a pair $(X,\Omega)$ where
  $X$ is a projective curve. Then there exists an snc model $\Yy$ and
  a dominant morphsm $f: \Yy \to \Xx$ such that $(\Yy, f^*
  \Omega')$ is an snc model of the pair $(X, \Omega)$
\end{lemma}

\begin{proof}
  Suffices to take $\Yy$ to be a log-resolution of $\Xx_0 \cup
  \div(\Omega)$. Indeed, $\div(f^* \Omega) = f^* \div(\Omega)$.
\end{proof}

Admitting a slight abuse of notation, given a $K$-varity $X$ and its
model $\Xx$ over $R$ we will denote by $X_s$ (resp. $\Xx_s$) the
fibres for $s$ close enough to 0 of the corresponding families over a
disc (resp. punctured disc). We will also denote $\omega_0$ the
$(1,1)$-form $\dfrac{i}{2}(\Omega_0 \wedge \bar\Omega_0)$ defined on
the smooth part of $\Xx_0$.

\begin{prop}
  \label{expansion}
  Let $(\Xx,\Omega')$ be an snc model of a pair $(X,\Omega)$.  Let
  $\sigma$ be an edge and $\{i,j\} \in \partial \sigma$. Assume that
  $\div(\Omega) \cup E_i \cup E_j$ is a divisor with snc support. Let
  $x_\sigma, y_\sigma \in \Oo_{\Xx_s,\sigma}$ be the local equations
  of $E_i, E_j$, respectively, $x_\sigma^{N_i} y_\sigma^{N_j} = t$ be
  the local equation of $E_i \cup E_j$ near $\sigma$, $a,b$ are
  integers, $a \geq 0, b \leq 0$ such that $aN_j + bN_i=(N_i,N_j)$ and
  put $z_\sigma = x_\sigma^a y_\sigma^{-b}$. Then for $s$ sufficiently
  close to 0
  $$
  \dfrac{i}{2}(\Omega'_s \wedge \bar\Omega'_s) =
  |s|^{2(\wt_\Omega(E_j) +\frac{b}{N_j}k -1)} |z_\sigma|^{2k}
  |u|^2 \dfrac{d|z_\sigma|}{|z_\sigma|} \wedge d\Arg z_\sigma
  $$
  where
  $z_\sigma = x_\sigma^a y_\sigma^{-b},
  k=\dfrac{(\wt_\Omega([E_i])-\wt_\Omega([E_j]))N_i N_j}{(N_i,N_j)}$ and
  $u \in \Oo^\times_{\Xx_0,\sigma}$.
\end{prop}

\begin{proof}
  Denote $c=\ord_{E_i}(\Omega'), d = \ord_{E_j}(\Omega')$. By a
  standard computation (see, for example, \cite[4.1.4]{mn-weight}) the
  form $\dfrac{dy}{x_\sigma^{N_i-1}y_\sigma^{N_j}}$ generates
  $\omega_{\Xx/R}$. Since
  $$
  \dfrac{dz_\sigma}{z_\sigma} = a\dfrac{dx_\sigma}{x_\sigma} -
  b\dfrac{dy_\sigma}{y_\sigma} 
  = (b - \dfrac{a N_j}{N_i})\dfrac{dy_\sigma}{y_\sigma} + \dfrac{a
    dt}{N_i t} 
  $$
  it follows that
  $$
  \div(\dfrac{dz_\sigma}{z_\sigma}) =
  \div(\dfrac{dy_\sigma}{y_\sigma}) = (N_i - 1) E_i + (N_j - 1)E_j
  $$
  Then
  $\Omega'_s = u x_\sigma^{c - N_i + 1} y_\sigma^{d - N_j + 1}
  \dfrac{dz_\sigma}{z_\sigma}$ for some unit 
  $u \in \Oo^\times_{\Xx_0,\sigma}$. Since
  $$
  \dfrac{dz_\sigma}{z_\sigma} \wedge \dfrac{d\bar z_\sigma}{\bar
    z_\sigma} = -2i \dfrac{d|z_\sigma|}{|z_\sigma|} \wedge d\Arg
  z_\sigma
  $$
  it is left to show that
  $$
  x_\sigma^{c - N_i + 1} y_\sigma^{d - N_j + 1} = s^{\wt_\Omega([E_j]) +\frac{b}{N_j}k-1} z_\sigma^k
  $$
  Note that by definition
  $\wt_{\Omega}([E_i]) = \dfrac{1+c}{N_i}, \wt_{\Omega}([E_j]) =
  \dfrac{1+d}{N_j}$ and
  $b\dfrac{N_j}{N_i} + a)k = (\wt_{\Omega}([E_i]) -
  \wt_{\Omega}([E_j]))N_i$. Expanding and simplyfying the rhs we get
  \begin{dmath*}
  x_\sigma^{\wt_{\Omega}([E_j])N_i+b\frac{N_i}{N_j}k- N_i}
  y_\sigma^{1+d + bk - N_j} x^{ak} y^{-bk} \\ =
  x_\sigma^{\wt_{\Omega}([E_j])N_i - N_i} 
  y_\sigma^{1+d-N_j} = x_\sigma^{1+c-N_i} y_\sigma^{1+d-N_j}
  \end{dmath*}
\end{proof}

\begin{cor}
  \label{omega-regular}
  Let $(\Xx,\Omega')$ be an snc model of a pair $(X,\Omega)$. Assume
  $\wt_\Omega([E_i]) < \wt_\Omega([E_j])$ for all $j \in \St(i)$. Then
  $\Omega/s^{\wt_\Omega([E_j])+b/N_j-1}|_{E_i}$ is regular and non-zero.
\end{cor}

We will need an snc model of $X$ that satisfies the following
technical assumption in order to describe the collapsed limit in
Theorem~\ref{limit-collapsed}:
\begin{assumption}
  \label{irred} 
  For all prime divisors $E_i, E_j \subset \Xx_0$ and for each
  $\sigma \in E_i \cap E_j$ there exists a neighbourhood $U$ of
  $\sigma$ such that $U_t$ is connected for $t$ close enough to 0.
\end{assumption}

\begin{lemma}
  For any pair $(X,\Omega)$ there exists a finite extension
  $K' \supset K$ such that the pair $(X \otimes K', \Omega \otimes K')$
  has an snc model that satisfies Assumption~\ref{irred}.
\end{lemma}

\begin{proof}
  By semi-stable reduction \cite{kkmd} there exists a finite extension
  $K' \supset K$ such that $X \otimes K'$ has an snc model $\Xx$ with
  all the irreducible components of $\Xx_0$ having multiplicity 1, so
  the Assumption~\ref{irred} is satisfied. Then notice that
  Assumption~\ref{irred} is stable under blow-ups and apply
  Lemma~\ref{make-snc}.
\end{proof}

\begin{prop}
  \label{charts}
  Let $(\Xx, \Omega')$ be an snc model of a pair $(X, \Omega)$ that
  satisfies Assumption~\ref{irred}. Then there exists a cover
  $U^\alpha \subset X(\C)$ of a neigbourhood of $\Xx_0$ that is
  indexed by vertices and edges of $\Delta_\Xx$ that satisfies the
  folowing properties:
  \begin{enumerate}[label=(\roman*)]
  \item\label{form} for all edges $\sigma$, the set $U^\sigma$ is
    defined in a neighbourhood of $\sigma$ by the inequalities
    $$
    C_\sigma |s|^{a/N_i} \leq |z_\sigma| \leq  D_\sigma |s|^{-b/N_j}
    $$
    for some constants $C_\sigma, D_\sigma > 0$, $U^\sigma_s$ is
    connected for all $s$ close enough to 0;
  \item\label{int} for any vertex $i$ and any edge $\sigma$,
    $U^i \cap U^\sigma \neq \emptyset$ if and only if
    $i \in \partial \sigma$;
  \end{enumerate}
\end{prop}

\begin{proof}
  Let $\epsilon$ be a number such that
  $x_\sigma^{N_i} y_\sigma^{N_j}=s$ are the equations of $\Xx_s$ in
  the total space of the degeneration for $|s| < \epsilon$ for all
  edges $\sigma \in \Delta_\Xx$.

  Observe that the inequalities from the property~\ref{form} can be
  rewritten as $|x_\sigma| \leq \dfrac{1}{C_\sigma^{N_j/(N_i,N_j)}}$
  and $|y_\sigma| \leq \dfrac{1}{D_\sigma^{N_i/(N_i,N_j)}}$, and pick
  the constants $C_\sigma, D_\sigma$ so that
  $U^\sigma \cap U^\tau \neq \emptyset$ for each vertex $i$ and each
  $\sigma, \tau \in \St(i)$. The property~\ref{form} follows from
  Assumption~\ref{irred}.

  The complement of the union of $U^\sigma$ over all edges $\sigma$ in
  the neighbourhood of the special fibre defined by the inequality
  $|s| < \epsilon$ consists of connected components $W^i$ that are in
  bijective correspondence with irreducible components
  $E_i \subset \Xx_0$. Define
  $U^i = W^i \cup \bigcup_{\sigma \in \St(i)} \partial U^\sigma$. Then
  the property~\ref{int} is satisfied by construction.
\end{proof}

\begin{lemma}
  \label{unit-bound}
  Let $u \in  \Oo^\times_{\Xx, \sigma}$ and let $W$ be a set
  defined by the inequalities
  $|s|^{\alpha} \leq |z_\sigma| \leq |s|^\beta$, $|s| < \eps$ in the
  neighbourhood of $\sigma$. If $\alpha \leq a/N_i$,
  $\beta \geq -b/N_j$ then $\sup_{x \in W_s} |u| = O(1)$, and if
  $\alpha < a/N_i$, $\beta > -b/N_j$ then
  $\sup_{x \in W_s} |u| = o(1)$, as $|s| \to 0$.
\end{lemma}

\begin{proof}
  The conclusion follows immediately after observing that the set
  $W_s$ is the intersection of the curve $x^{N_i} y^{N_j} = s$ and the
  rectangle $|x| \leq |s|^{a/N_i - \alpha}$, $|y| \leq |s|^{\beta + b/N_j}$.
\end{proof}

\subsection{Asymptotic distance estimates}

\begin{lemma}
  \label{annuli}
  Let $(\Xx, \Omega')$ be an snc model of the pair
  $(X,\Omega)$. Consider the fibres $X_s$ with the K\"alher metric
  $\omega_s = \dfrac{i}{2}(\Omega_s \wedge \bar\Omega_s)$. Pick
  functions $z_\sigma \in \Oo_{\Xx, \sigma}$ as in
  Proposition~\ref{expansion} and assume that
  $$
  \Omega' = c_\sigma z_\sigma^\alpha s^\beta (1 + u)dz_\sigma
  $$
  for some $\alpha, \beta \in \Z$, and
  $u \in \m \subset \Oo_{\Xx_s, \sigma}$.
  
  Let $\{\eta_s,\eta'_s\}$ be a collection of points in $U^\sigma_s$
  and let $\gamma_s \subset U^\sigma_s$ be a shortest path between the
  points $\eta_s$ and $\eta'_s$. Then the length of $\gamma_s$ as
  $s \to 0$ is
  $$
  |c_\sigma|\cdot |\ln
  |z_\sigma(\eta'_s)| - \ln |z_\sigma(\eta_s)||\cdot |s|^\beta (1 + o(1))
  $$
  if $\alpha=-1$ and 
  $$
  (|z_\sigma(\eta'_s)|^{\alpha+1}-|z_\sigma(\eta_s)|^{\alpha+1})\dfrac{|s|^\beta}{\alpha+1}(1
  + o(1))
  $$
  otherwise.
\end{lemma}

\begin{proof}
  The Riemannian metric tensor on in polar coondinates $\rho=|z_\sigma|,
  \theta=\Arg z_\sigma$ is given by the expression
  \begin{dmath*}
    g_s = |c_\sigma|^2 |z_\sigma|^{2\alpha}|s|^{\beta}|1+u|^2\, (d\rho
    \otimes d\rho + \rho^2 d\theta \otimes d\theta)
  \end{dmath*}
  in the neighbourhood of $\sigma$.  Let
  $\gamma'_s:[0,1] \to U^\sigma_s$ be the path given by
  $$
  \gamma'_s(\tau) = ((1-\tau) |z_\sigma(\eta_s)| + \tau 
  |z_\sigma(\eta'_s)|) \exp( i \Arg(\eta_s))
  $$
  and $\gamma''_s: [0,1] \to U^\sigma_s$ be defined by
  $$
  \gamma''_s(\tau) = |z_\sigma(\eta'_s)| \exp(i (\tau \Arg \eta_s + (1 -
    \tau) \Arg \eta'_s))
  $$
  Clearly,
  $$
  \dfrac{d}{d\tau}\gamma'(\tau) = (|z_\sigma(\eta')| -
  |z_\sigma(\eta_s)|) d\rho \qquad \dfrac{d}{d\tau}\gamma''(\tau) =
  (\Arg \eta'_s - \Arg \eta_s)  d \theta
  $$  
  
  Assume for definiteness that
  $|z_\sigma(\eta'_s)| > |z_\sigma(\eta_s)|$. Then
  $$
    L(\gamma'_s) = \displaystyle\int_0^1
    \sqrt{g_s(\dot{\gamma'}(\tau), \dot{\gamma'}(\tau)} d\tau =
    \displaystyle\int_{|z_\sigma(\eta_s)|}^{|z_\sigma(\eta'_s)|}|c_\sigma|
    \rho^\alpha|s|^\beta|1 + u(\rho e^{i\Arg \eta_s},s)| d\rho
  $$
  Denote
  $I_{\epsilon,s} = [|s|^{a/N_i+\epsilon}, |s|^{-b/N_j-\epsilon}]$ and
  let $H_s=[|z_\sigma(\eta_s)|, |z_\sigma(\eta'_s)|]$; denote
  $A_{\epsilon,s}, B_{\epsilon,s}$ the endpoints of the interval
  $I_{\epsilon, s} \cap H_{s}$.  The latter integral can be
  represented, for $\epsilon > 0$, as the sum of two integrals
  \begin{dmath*}
    \displaystyle\int_{H_s} |c_\sigma| \rho^\alpha|s|^\beta|1 +
    u(\rho e^{i\Arg \eta_s}, s)| d\rho = |c_\sigma| |s|^\beta \left(
      \displaystyle\int_{H_s\wo I_{\epsilon, s}} \rho^\alpha |1 +
      u(\rho e^{i\Arg \eta_s}, s)| d\rho + \\ + \displaystyle\int_{I_{\epsilon, s}}
      \rho^\alpha |1 + u(\rho e^{i\Arg \eta_s}, s)| d\rho \right)
  \end{dmath*}
  Let us first consider the case $\alpha=-1$.  By
  Lemma~\ref{unit-bound}
  \begin{dmath*}
    L(\gamma'_s)  = |c_\sigma| |s|^\beta \sup_{\epsilon > 0}\, (\ln
    B_{\epsilon,s}/A_{\epsilon,s})(1+o(1)) + \\ + C \max
    \{A_{\epsilon,s} -    
    |z_\epsilon(\eta_s)|,0\} + C \max \{|z_\epsilon(\eta'_s)| -
    B_{\epsilon,s},0\}) = |c_\sigma| |s|^\beta \ln
    |z_\sigma(\eta'_s)|/|z_\sigma(\eta_s)|(1+o(1))
  \end{dmath*}
  Similarly, one derives for $\alpha \neq -1$,
  $$
  L(\gamma'_s) =
  \dfrac{|s|^{\beta}}{\alpha+1}(|z_\sigma(\eta'_s)|^{\alpha+1} -
  |z_\sigma(\eta_s)|^{\alpha+1})(1+o(1))
  $$    
  Clearly, $L(\gamma_s) \leq L(\gamma'_s) + L(\gamma''_s)$ and
  $L(\gamma''_s) = O(|s|^{\beta}|z_\sigma(\eta'_s)|)$.

  On the other hand, denoting the polar coordinates of $\gamma_s$ by
  $\gamma_{s,\rho}, \gamma_{s,\theta}$ we have
  \begin{dmath*}
  L(\gamma_s) = \int_0^1
  \sqrt{g_s(\dot{\gamma}(\tau),\dot{\gamma}(\tau)} d\tau =
  \int_0^1 \sqrt{|c_\sigma|^2 (\gamma_{s,\rho})^{2\alpha} |s|^{2\beta}
    |1+u(\gamma_s)|^2 ((\dot{\gamma}_{s,\rho})^2 +
    \gamma_{s,\rho}^2(\dot{\gamma}_{s,\theta})^2 } d\tau 
  \geq \int_0^1 |c_\sigma| (\gamma_{s,\rho})^{\alpha} |s|^{\beta}
    |1+u(\gamma_s)| (\dot{\gamma}_{s,\rho}) d\tau 
  \end{dmath*}
  The last expression has the same asymptotics as $L(\gamma'_s)$ and we
  conclude.
\end{proof}

\begin{lemma}
  \label{outside-annuli}
  Let $(\Xx, \Omega')$ be an snc model of a pair $(X,\Omega)$ and let
  $E_i$ be an irreducible component of $\Xx_0$. Let
  $a_s, b_s \in \Xx_s, |s| < \eps $ be collections of points such that
  $\lim_{s \to 0} a_s = a_0, \lim_{s \to 0} b_s = b_0$ for some
  $a_0, b_0 \in E_i \subset \Xx_s$. Assume that $\wt_\Omega([E_i]) =
  1$, $\wt_\Omega([E_j]) > 1$ for all
  $[E_j] \in \St([E_i])$
  and that $a_0, b_0 \notin (\Xx_0)_{sing}$. Let $\gamma_s \subset U^i$
  be a shortest path for the metric
  $\omega_s=\dfrac{i}{2}(\Omega_s \wedge \bar\Omega_s)$ on $X_s$ and
  that connects $a_s$ and $b_s$. Then
  $$
  \lim_{s\to 0} l(\gamma_s) = l(\gamma_0)
  $$
  and the limit is finite.
\end{lemma}

\begin{proof}
    As was observed before in the proof of Lemma~\ref{expansion} the
  form $\Omega'_0|_{E_i}$ can be written down in a Zariski
  neighbourhood of $\sigma \in E_i \cap E_j$ as
  $$
  \Omega'_0|_{E_i} = u x_\sigma^{(\wt_\Omega([E_i]) - 1)N_i}
  y_\sigma^{(\wt_\Omega([E_j]) - 1)N_j} \frac{dy_\sigma}{y_\sigma}
  $$
  where $u \in \Oo^\times_{\Xx, \sigma}$ and $x_\sigma, y_\sigma$ are
  local equations of $E_i, E_j$ respectively. We may assume that $E_i$
  is covered by such neighbourhdoods, passing to a model where
  finitely many points of $E_i$ are blown up if it is not the case;
  notice that in this case the assumptions about the weight remain
  true.  Since $(\wt_\Omega([E_j]) - 1)N_j > 0$ and integral,
  $\Omega'_0|_{E_i}$ is regular and non-zero. The statement then
  follows from the continuity of the form
  $\dfrac{i}{2} \Omega' \wedge \bar\Omega'$.
\end{proof}

\begin{lemma}
  \label{diam-order}
  Let $(\Xx, \Omega')$ be an snc model of a pair $(X,\Omega)$ that
  satsfies Assumption~\ref{irred}.  Let
  $k=\min_{x \in \Delta_\Xx} \wt_\Omega(x)$ and let
  $\Sigma \subset \Delta(\Xx_s)$ be the set of faces of $\Delta_\Xx$
  where $\wt_\Omega \equiv k$. Let $\Omega''=\Omega'/s^{k-1}$ and let
  $\omega''_s=\dfrac{i}{2}(\Omega''_s \wedge \bar\Omega''_s)$. The
  asymptotics of the diameter of $X_s$ as $s \to 0$ with respect to
  the K\"ahler metric
  $\omega_s=\dfrac{i}{2}(\Omega_s \wedge \bar\Omega_s)$ is
  $$
  \diam X_s = (\diam (X_0,\omega''_0) |s|^{k-1} (1 + o(1))
  $$  
  if $\dim \Sigma=0$, and 
  $$
  \diam X_s = c \ln |s| |s|^{k-1} (1+o(1))
  $$
  if $\dim \Sigma=1$, for some constant $c$.
\end{lemma}

\begin{proof}
  Consider the cover constructed in Proposition~\ref{charts}.
   
  Since $\wt_{\Omega'/t^m}=\wt_{\Omega'}-m$ for any integer $m$, by
  Lemma~\ref{expansion} 
  $$
  \diam U^i_s = (\diam (U^i_0,) |s|^{\wt_\Omega([E_i])-1}(1 + o(1))
  $$
  for any vertex $[E_i] \in \Delta_\Xx$.

  By the same consideration, after applying Lemmas~\ref{annuli} and
  \ref{outside-annuli} we get
  $$
  \diam U^\sigma_s = d_\sigma\ln|s|\,|s|^{\wt_\Omega([E_i])-1}(1 + o(1))
  $$
  when $\wt_\Omega([E_i])=\wt_\Omega([E_j])$ for $[E_i],[E_j] \in \partial
  \sigma$,
  $$
  \diam U^\sigma_s = d_\sigma|s|^{\wt_\Omega([E_i])-1}(1 + o(1))
  $$
  when $\wt_\Omega([E_i]) < \wt_\Omega([E_j])$, for some positive real
  constants $d_\sigma$.

  It follows that $\diam((U^\beta)_s) = o(\diam (U^\alpha)_s)$ for any
  $\alpha \in \Sigma$ and $\beta \notin \Sigma$. Now notice that
  $\diam U^\sigma$ has the same asymptotics as $|s| \to 0$ for all
  edges $\sigma \in \Sigma$ if $\dim \Sigma=1$, and that $\diam U^i$
  have the same asymptotics for all $i \in \Sigma$ if $\dim
  \Sigma=0$. Therefore, $\diam (X_s, \omega''_s)$
  in the first case tends to the diameter of $(X_0,\omega''_0)$, and
  $\diam (X_s, \omega''_s)$ tends to some positive real
  number (the diameter of the graph $\Delta_\Xx$ with lengths of edges
  adjusted) in the second case. The conclusion of the Lemma follows.
\end{proof}

To describe the Gromov-Hausdorff limit of a family of curves
$(X_s, \tilde{\omega}_s)$ we will distinguish two cases:
\emph{collapsed} limit, when the diameter of $(X_s, \omega_s)$ is of
order $(\ln |t|\,|t|^k)(1+o(1))$ for some $k$, and \emph{non-collapsed}
limit otherwise.

\subsection{Shape of the limit}
\label{shape}

If $X$ and $Y$ are two (pseudo)metric spaces and
$R \subset X \times Y$ is a relation one defines the \emph{distortion}
of $R$ to be
$$
\dis R = \sup\limits_{(x,y), (x',y') \in R} |d_X(x,y) - d_Y(x',y')|
$$
One easily observes that if both projections of $R$ on $X$ and $Y$ are
surjective then $d_{GH}(X,Y) \leq \dis R / 2$, and conversely, for any
metric spaces $X, Y$ such that $d_{GH}(X,Y) \leq \epsilon$ the
relation $R_\eps = \{ (x,y) \in X \times Y \mid d(x,y) < \eps \}$,
where $d$ is the metric on $X \sqcup Y$ that realizes the bound,
satisfies $\dis R_{\eps} < 2\epsilon$ (\cite[Theorem~7.3.25]{bbi}).

Recall that a metric is called an \emph{inner metric} if the distance
beween two points is defined as an infimum of a length functional on
some class of admissible paths (see Section~2 of \cite{bbi} for the
detailed definition).

\begin{lemma}
  \label{patching}
  Let $(X, d)$ and $(Y, d')$ be two pseudometric spaces, let
  $\cup_{i=1}^n U_i = X$ and $\cup_{j=1}^n V_j = Y$ be two coverings
  by path-connected sets, and let $R \subset X \times Y$ be a
  relation.  Assume that
  \begin{enumerate}[label=(\roman*)]
  \item the 1-nerves of $\{ U_i \}$ and $\{ V_j \}$ are isomorphic,
    that is, for all $i\neq j$ the connected components $U_\sigma$ of
    $U_i \cap U_j$ are in bijective correspondence with connected
    components $V_\sigma$ of $V_i \cap V_j$;
   \item for all $i,j$ (including $i=j$), for all connected components
     $U_\sigma \subset U_i \cap U_j $ the relation
     $R \cap (U_\sigma \times V_\sigma)$ projects surjectively on
     $U_\sigma, V_\sigma$;
   \end{enumerate}
   Then there exists a number $N$, depending only on the nerve of
   $\{U_i\}$ and $\{V_i\}$ such that 
   $$
   \dis R \leq \max_{i_0, \ldots, i_N} \sum_{k=1}^N \dis R \cap (U_{i_k} \times V_{j_k})
   $$
   where the maximum is taken over such sequences $\{i_k\}$ that
   $U_{i_k} \cap U_{i_{k+1}} \neq \emptyset$ for all $k$.
 \end{lemma}
 
 \begin{proof}
   Let $N$ be a number such that any path in $X$ (or, equivalently,
   $Y$) passes consecutevely through a sequence of elements of the
   cover $U_{i_0}, \ldots, U_{i_L}$, $L \leq N$, so that any finite
   subsequence starting and ending with the same element occurs at
   most once.
   
   Let $x_0, x_L \in X, x_0 \in U_{i_0}, x_L \in U_{i_L}$ and
   subdivide the shortest path between $x_0$ and $x_L$ by adding
   points $x_k$ so that $x_k,x_{k+1} \in U_{i_k}$ for some
   sequence $\{i_k\}$. Pick $y_0, \ldots, y_L$ so that
   $y_k, y_{k+1} \in V_{i_k}, (x_i, y_i) \in R$. Then
   \begin{dmath*}
   d'(y_0, y_L) \leq \sum_{k=0}^L d'(y_k,y_{k+1}) \leq
   \sum_{i=0}^L d(x_k,x_{k+1}) + \dis R \cap (U_{i_k} \times V_{i_k})
   \leq 
   d(x_0, y_L) + \sum_{k=0}^L \dis R \cap (U_{i_k} \times V_{i_k})
   \end{dmath*}
   By the symmetric argument we obtain that
   $$
   d(x_0, x_L) \leq d'(y_0, y_L) + \sum_{k=0}^L \dis R \cap (U_{i_k}
   \times V_{i_k})
   $$

   Since $L \leq N$, we can conclude.
 \end{proof}

\begin{thm}[Collapsed limit]
  \label{limit-collapsed}
  Let $(\Xx,\Omega')$ be an snc model of $(X,\Omega)$ that satisfies
  Assumption~\ref{irred}. Let
  $k=\min_{x \in \Delta_\Xx} \wt_\Omega(x)$ and let
  $\Sigma \subset \Delta(\Xx_s)$ be the union of vertices and edges of
  $\Delta_\Xx$ where $\wt_\Omega \equiv k$.  Assume that
  $\dim \Sigma=1$.  Let $x \sim y$ for $x,y \in \Delta_\Xx$ if and
  only if there exists a path $\gamma: [0,1] \to \Delta_\Xx$ joining
  $x$ and $y$ such that $|\gamma^{-1}(\Sigma)| < \infty$.  The
  Gromov-Hausdorff limit of $(X_s, \tilde{\omega}_s)$ as $s \to 0$ is
  isometric to $\Delta_\Xx/\sim$ endowed with the metric that
  stretches each edge $\sigma$ by the factor $|c_\sigma|$ and
  renormalized so that $\diam (\Delta_\Xx/\sim) = 1$.
\end{thm}

\begin{proof}
  We adopt the notation for local coordinates from
  Proposition~\ref{charts}. By Lemma~\ref{diam-order},
  $\diam X_s = c \ln |s| |s|^{k-1} (1 + o(1))$ for some constant
  $c$. We identify each edge $\sigma$ with an interval
  $[-\dfrac{b|c_\sigma|}{N_jc}, \dfrac{a|c_\sigma|}{N_ic}]$. Define
  the map $f_s: X_s \to \Delta_\Xx$ as follows:
  $$
  f_s(x,s) = \left \{
    \begin{array}{ll}
      \dfrac{|c_\sigma|}{c}\, \Xi_\sigma\left(\dfrac{\ln |z_\sigma|}{\ln |s|}\right) \in [\sigma], & 
        \textrm{ if } x \in U^\sigma\\
       {[ E_i ]}  & \textrm{ if } x \in U^i\\
    \end{array}
  \right.
  $$
  where $\Xi_\sigma: [-b/N_j - \ln D_\sigma, a/N_i - \ln C_\sigma] \to
  [-b/N_j, a/N_i]$ is the linear bijection.
  Let $R_s$ be the graph of $f_s$. It follows from
  Assumption~\ref{irred} and the fact that $f_s|_{U^\sigma_s}$ is
  surjective onto $[\sigma]$ that $R_s$ satisfies the conditions of
  Lemma~\ref{patching}.

  Consider the metric
  $\tilde{\omega}=\dfrac{\omega_s}{\diam X_s^2}$ on $X$.  It
  follows from Lemma~\ref{annuli} that
  $\dis  R_s \cap (U^\sigma \times [\sigma]) \to 0$ as $s \to 0$
  for $\sigma \in \Sigma$.  Since by Lemma~\ref{diam-order}
  $\diam U^\sigma_s = o(\diam X_s)$ for any
  $[\sigma] \not\subseteq \Sigma$,
  $\dis  R_s \cap (U^\sigma \times [\sigma]) \to 0$ for such
  $[\sigma]$. Since by Lemma~\ref{outside-annuli}
  $\diam U^i_s \to 0$, $R_s \cap (U^i \times [E_i]) \to 0$
  as $s \to 0$.

  Therefore, by Lemma~\ref{patching}, $\dis R_s \to 0$ as $s \to 0$.
  If $q: \Delta_\Xx \to \Delta_\Xx/\sim$ is the projection on the
  quotient then clearly $\dis q=0$, since $\Delta/\sim$ is the metric
  space associated to the pseudometric space $\Delta$. It follows that
  if $\bar R_s$ is the graph of $q \circ f_s$ then $\dis R_s \to 0$ as
  $s \to 0$. It follows that $X_s$ converges in the sence of
  Gromov-Hausdorff to $\Delta_\Xx/\sim$.
\end{proof}

\begin{rem}
  The metric graph $\Delta_\Xx/\sim$ does not depend on the choice of
  a model $\Xx$. Indeed, observe that if $\Yy$ is a model that
  dominates $\Xx$ then $\Delta_\Yy$ contracts onto $\Delta_\Xx$ and it
  follows from Proposition~4.3.4 \cite{mn-weight} that
  $\Sigma_\Yy=\Sigma_\Xx$, then use the fact that any two models are
  related by a series of blow-ups and blow-downs of points in
  the special fibre.
\end{rem}

As was observed in \cite[Lemma~3.4.5]{bn-weight}, if the genus of
$E_i$ 0 then $[E_i] \in \Delta_\Xx$ belongs to $\Sigma$ if and only if
some adjacent edge belongs to $\Sigma$. In the case of non-collasped
limit the components $E_i$ of the central fibre such that
$[E_i] \in \Sigma$ can thus be regarded as surfaces endowed with a
flat metric, since by Lemma~\ref{outside-annuli} the restriction of
$\Omega'/t^{k-1}$ to each such $E_i$ is a regular 1-form.

\begin{thm}[Non-collapsed limit]
  \label{limit-non-collapsed}
  Let $(\Xx,\Omega')$ be an snc model of $(X,\Omega)$ as in the
  previous theorem and assume that $\dim \Sigma=0$.  Let $\sim$ be the
  smallest equivalence relation on $X_0$ containing the relation
  defined as follows:
  \begin{enumerate}
  \item $\sigma \sim \tau$ for $\sigma \subset E_i$,
    $\tau \subset E_j$ if $[E_i],[E_j] \in \Sigma$ and there exists a
    path $\gamma: [0,1] \to \Delta_\Xx$ such that
    $\gamma(0) = [E_i], \gamma(1) = [E_j]$, initial segment of
    $\gamma$ passes through $[\sigma]$, and final segment of $\gamma$
    passes through $[\tau]$,
  \item $x \sim \sigma$ if $i \notin \Sigma$, $x \in E_i$
    $\partial [\sigma] \cap \Sigma \neq \emptyset$, and there exists a
    path $\gamma: [0,1] \to \Delta_\Xx$ such that $\gamma(0) = [E_i]$,
    and $[\sigma] \subset \gamma([0,1])$.
  \end{enumerate}
  Then the limit of $(X_s, \tilde{\omega}_s)$ as $s \to 0$ is
  $X_0 / \sim\,=(\cup_{[E_i] \in \Sigma} E_i) / \sim$ with the metric
  renormalized so that the diameter of the space is 1.
\end{thm}

\begin{proof}
  The central fibre $\Xx_0$ is a deformation retract of the total
  space $X$, (see \cite{clemens} for an explicit construction of a
  retraction). Let $r: X \to \Xx_0$ be a retraction, and denote $r_s$
  its restriction to $X_s$. Let $R_s$ be the graph of $r_s$.

  Then by Lemmas~\ref{diam-order} and \ref{outside-annuli},
  $\dis R_n \to 0$ as $s \to 0$ if we consider $X_0$ as a
  pseudo-metric space with the metric given by the form
  $\Omega'/t^{k-1}$. The relation $\sim$ described in the statement of
  the Theorem identifies points at distance 0 between each other. It
  then follows that $X_0 / \sim$ is the Gromov-Hausdorff limit of
  $X_s$ as $s \to 0$.
\end{proof}

 

We will now study for the purpose of illustration of
Theorem~\ref{limit-collapsed} the possible shapes of the
Gromov-Hausdorff limits it describes.

We will use the description of the graph Laplacian of the weight
function due to Baker and Nicaise \cite{bn-weight}, which we quickly
recall. By a \emph{weighted graph} we understand a metrized graph
$\Gamma$ with the set of vertices $V(\Gamma)$ and with infinite edges
allowed, and a pair of functions $N,g : V(\Gamma) \to \Z$. Given an
snc model $(\Xx, \Omega)$ of $(X, \Omega)$, one associates a weighted
graph as follows: take $\Delta_\Xx$ and attach infinite edges at the
vertices which correspond to components having non-trivial
intersection with $\div(\Omega)$. A \emph{divisor} on $\Gamma$ is a
formal combination of the vertices of $\Gamma$. Let $f: \Gamma \to \R$
be a function that is affine on every edge of $\Gamma$, then the
\emph{Laplacian of $f$} is the divisor
$\Delta(f) = \sum_{i \in V(\Gamma)} a_i v_i$ where $a_i$ is the sum of
outward slopes of $f$ at $v_i$. The \emph{canonical divisor of
  $\Gamma$} is the divisor
$$
K_\Gamma = \sum_v N_v( \val(v) + 2g(v) - 2) v
$$
where $\val(v)$ is the valency of the vertex $v$. By
\cite[Theorem~3.2.3]{bn-weight} each infinite edge running from a
vertex $v$ towards a zero $x \in X$ of a differential form $\Omega$
has an outgoing slope $N_v(1 + \deg_x(\Omega))$ and
$\Delta(\wt_\Omega)=K_\Gamma$.

\begin{prop}
  \label{collapsed-any-genus}
  For any $k \geq 1$ the wedge sum of $k$ circles can occur as a limit
  (in the sense of Theorem~\ref{limit-collapsed}) of a family of
  curves of genus $2k+1$. This family admits an snc model such that
  all the components of the central fibre are rational.
\end{prop}

\begin{proof}
  For the case $k=1$ take an elliptic curve over $\Cmer$ of bad
  reduction and any regular form $\Omega$.

  For $k>1$, we will construct a weighted graph $\Gamma^+$ with a
  subgraph $\Gamma$ of non-infinite edges of $\Gamma^+$, and a weight
  function $w: \Gamma^+ \to \R$ affine on each edge, we will then
  apply \cite[Theorem~6.3]{muw} to get a pair $(X,\Omega)$ that gives
  rise to the Gromov-Hausdorff limit of the desired genus.

  Let $\Gamma$ be a chain of $g=2k-1$ cycles, $C_1, \ldots, C_{2k-1}$
  joined consecutively, each cycle $C_i$ consisting of two edges
  connecting vertices $c_i$ and $c_{i+1}$. Attach two infinite edges
  (corresponding to the zeroes of the differential form) to both edges
  of each $C_{2i}, 1 \leq i \leq k$, subdividing them in the points of
  attachment $a_{2i-1}, a_{2i}, b_{2i-1}, b_{2i}, 1 \leq i \leq k$;
  call the resulting graph $\Gamma^+$. Define the function
  $w: V(\Gamma^+) \to \Z$ as follows:
  $$
  w(c_i) = 0 \ \ w(a_i)=w(b_i)=1
  $$
  Extend it in the affine fashion to the edges of $\Gamma^+$, with the
  outgoing slope 2 on the infinite edges. Put $N=0, g=0$ on all
  vertices, and let all edges be of length 1. One checks that
  $\Delta(\wt_\Omega)=K_{\Gamma^+}$. Clearly, the minimality locus of
  $w$ is the union of odd cycles $C_{2i-1}$, and the quotient by the
  equivalence relation described in Theorem~\ref{limit-collapsed} is of
  genus $k$.

  In order to apply \cite[Theorem~6.3]{muw} we need to find a
  piecewise affine function $f: \Gamma \to \R$ affine on the edges
  such that the tropical divisor $\sum a_{i} + b_i$ can be presented
  as as $K_\Gamma+\Delta(f)$. We claim that such $f$ can be taken to
  be $-w|_\Gamma$. We can then check that the extended level graph
  $\Gamma^+(f)$ does not have inconvenient vertices (Definition~6.2,
  \emph{loc.cit.})  because for all vertices one of the outgoing
  slopes is 0. Therefore the only condition of Theorem~6.3,
  \emph{loc.cit.} that is left to check is about the edges, and it is
  fulfilled, since all horizontal edges belong to a horizontal simple
  cycle. By Theorem~6.3, \emph{loc.cit.}  there exists a variety $X$
  over $\Cmer$ and a 1-form $\Omega$ on $X$ and an snc model
  $(\Xx, \Omega')$ of the pair $(X,\Omega)$ such that $\Delta_\Xx$ is
  isomorphic to $\Gamma$ and $\wt_\Omega=w$ on $\Gamma$.
\end{proof}

Let us conclude with two questions.
\vspace{1ex}\\
\noindent \textbf{Question 1}: can
Proposition~\ref{collapsed-any-genus} be proved (perhaps with some
additional conditions) for a fixed genus of the elements of the family
$\Xx$? for a given fixed $\Xx$, constructing appropriate $\Omega$?

\noindent \textbf{Question 2}: can one characterise non-collapsed limits
starting from the description of Theorem~\ref{limit-non-collapsed}?

\bibliography{flat}

\vspace{3ex}

\noindent {\sc Department of Mathematics\\
  KU Leuven\\
  Celestijnenlaan 200B \\
  B-3001 Leuven (Heverlee)\\
  Belgium\\}
{\tt dmitry.sustretov@kuleuven.be\\}

\end{document}